\documentclass[]{gtart}
\usepackage{amsmath,amsthm,amscd,amssymb,verbatim}
\usepackage{psfrag,graphicx,epsf} 
\usepackage{latexsym}

  \newtheorem{thm}{Theorem}[section]
  \newtheorem{lem}[thm]{Lemma}
  
  \newtheorem{cor}[thm]{Corollary}
  \newtheorem{prop}[thm]{Proposition}

  \newtheorem{question}[thm]{Question}
  \newtheorem{conj}[thm]{Conjecture}

  \theoremstyle{definition}
  \newtheorem{dfn}[thm]{Definition}

  \def\l{\ell}
  \def\lk{\ell k}

\newcommand{\be}{\begin{equation}}
\newcommand{\ee}{\end{equation}}
\newcommand{\bea}{\begin{eqnarray}}
\newcommand{\eea}{\end{eqnarray}}
\newcommand{\bmini}{\footnotesize\begin{center}\begin{minipage}{5.5in}}
\newcommand{\emini}{\end{minipage}\end{center}\normalsize}
\newcommand{\pf}{{\em Proof: }}

\newcommand{\real}{\mathbb{R}}

\newcommand{\eg}{{\em e.g.}}

\newcommand{\eps}{\epsilon}
\newcommand{\inv}{^{-1}}

\newcommand{\id}{{\mbox{{\sc Id}}}}
\newcommand{\T}{{\mathcal T}}
\newcommand{\U}{{\mathcal U}}
\newcommand{\CC}{{\mathcal C}}
\newcommand{\df}[1]{{\bf #1}}
\newcommand{\ith}{i^{th}}
\newcommand{\braid}{{b_{(n_1,\cdots,n_k)}}}
\newcommand{\word}[1]{\left(#1\right)}
\newcommand{\x}{\sf x}
\newcommand{\y}{\sf y}
\newcommand{\z}{\sf z}

\begin{document}

\title{Flowlines transverse to fibred knots and links}

\author{Robert Ghrist}
\address{Department of Mathematics, University of Illinois, Urbana IL,
  61801 USA}
\email{ghrist@math.uiuc.edu}
\thanks{RG supported in part by NSF Grant \# DMS-0134408.}

\author{Eiko Kin}
\address{Department of Mathematics, Kyoto University, 
Oiwake-cho Kitashirakawa Sakyo-ku Kyoto-shi,
Kyoto 606-8502 Japan}
\email{kin@kusum.kyoto-u.ac.jp}
\thanks{EK supported in part by JSPS Research Fellowships 
for Young Scientists}

\begin{abstract}
Let $K$ be a knot or link in $S^3$ which is fibred --- the complement
fibres over $S^1$ with fibres spanning surfaces. We focus on 
those fibred knots and links which have the following property: every vector
field transverse to the fibres possesses closed flow lines of all 
possible knot and link types in $S^3$. Our main result is that a large
class of fibred knots and links has this property, including all fibred 
non-torus 2-bridge knots. In general, sufficient conditions include a 
pseudo-Anosov type monodromy map and a sufficiently high degree of symmetry. 
\end{abstract}


\maketitlepage

%

\section{Introduction}
\label{sec_intro}

We consider fibrations of knot or links complements in $S^3$ from the point
of view of transverse vector fields. 

Periodic flowlines of a vector field on 
$S^3$ possess knotting and linking information which is intimately 
related to the dynamical properties of the flow. Examples abound 
of classes of systems for which simple dynamics
implicate simple knotted orbits (those whose complement cannot
possess a hyperbolic structure) and vice versa: see 
\cite{Wada} [Morse-Smale flows], \cite{Fra} [Smale flows], 
\cite{CMAN,FN} [integrable Hamiltonian flows], \cite{BW:1} 
[the Lorenz system], and \cite{EG} [flows tangent to plane fields]. 
In this paper, we demonstrate that complicated fibrations 
typically force complicated knot and link types in any transverse 
vector field  --- indeed, as complicated as can be imagined. 

Consider a knot or link $K$ in $S^3$ which is 
fibred; that is, there is a fibration $\pi:S^3-K\to S^1$ with fibre
a spanning surface $\Sigma$ for $K$. Choose any vector field 
$X$ transverse to the fibres of $\pi$ (in particular, $X$ must be
nonvanishing). Birman and Williams originally asked the question, 
{\em ``Which knot types are forced to exist as periodic orbits of $X$?''} 
The question splits into cases based on Thurston's classification
theorem for surface maps applied to the monodromy $\Phi:\Sigma\to\Sigma$
of $\pi$. If the complement of $K$ is geometrically ``simple''
(\eg, the unknot, torus knots, iterated torus knots), then only a 
finite set of knot types are forced to exist for all $X$ (again, typically 
unknots, torus knots, and iterations). On the other hand, if the 
complement of $K$ is ``hyperbolic'' (in the sense that the 
monodromy map is of pseudo-Anosov type), then there are always 
an infinite number of distinct knot types as periodic orbits in any
vector field $X$ transverse to $\pi$. 

The paper \cite{BW:2} carefully constructed a \df{template} (a branched 
surface with semiflow) in the 
complement of the figure eight knot which captures all knot and link 
types forced by a transverse vector field. It was there noted that, although 
many knot types were present, others were seemingly impossible to 
locate. To the contrary, the paper \cite{G} showed that this template
is \df{universal} --- every 
knot and link type can be found on this branched surface. It follows
that any vector field transverse to a fibration of the figure eight knot
complement must possess all knot and link types as orbits.
In this paper, we explore to what extent this property holds for other
fibred knots and links. 

\begin{dfn}
\label{def_UF}
A fibred knot/link $K\subset S^3$ is said to be \df{universally fibred} if
each vector field $X$ transverse to the fibration possesses closed 
orbits of all knot and link types.
\end{dfn}

An obvious prerequisite for being universally fibred is that the monodromy 
be sufficiently complex: pseudo-Anosov type. 
For in the case of periodic type, 
the forced links may have only a finite number of knot types
represented (by suspending the periodic monodromy). 

\begin{figure}[hbt]
\begin{center}
\includegraphics[angle=0,width=4.0in]{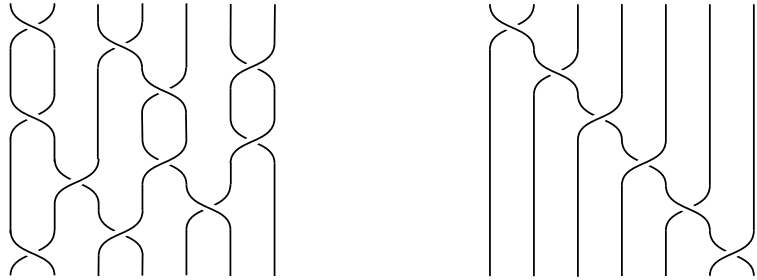}
\caption{[left] an example of a braid on seven strands; 
[right] an example of a braid represented by Eqn.~(\ref{eq_Main}), 
$b = \sigma_1\sigma_2\sigma_3^{-1}\sigma_4^{-1}\sigma_5\sigma_6^{-1}$} 
\label{fig_braid}
\end{center}
\end{figure}

Our results are phrased in the language of braids: see
\S\ref{sec_back} for details and Fig.~\ref{fig_braid}
for examples of braids. We summarize our results in a simplified 
form; see \S\ref{sec_maintheorems} for stronger results.

\vspace{0.1in}
\noindent
{\bf Main Theorem:} 
Let $b$ be a braid of the form
\begin{equation}
\label{eq_Main}
b = \sigma_1^{\eps_1}\sigma_2^{\eps_2}\cdots\sigma_N^{\eps_N},
\end{equation}
where $\eps_i=+1$ or $-1$ for each $i$ (see Fig.~\ref{fig_braid}[right]) and 
not all the $\eps_i$ are the same.  Then for every $m>1$, the 
closure of $b^m$ is a fibred knot/link in $S^3$ 
of pseudo-Anosov type. Furthermore:
\begin{enumerate}
\item[(a)] The closure of $b^m$ has the universal
fibration property for all but finitely many $m$.
\item[(b)] The closure of $b^2$ always has the universal 
fibration property.
\item[(c)] Whenever $b$ contains the subword 
$\sigma_{i-1}^{-1}\sigma_i\sigma_{i+1}^{-1}\sigma_{i+2}$ 
(or $\sigma_{i-2}\sigma_{i-1}^{-1}\sigma_i\sigma_{i+1}^{-1}$), 
then the closure of
$b^m$ is universally fibred for all $m>1$. 
\end{enumerate}
\vspace{0.1in}

 From result (b) above, we obtain a complete classification in 
the case of 2-bridge knots:

\vspace{0.1in}
\noindent
{\bf Corollary:} Every fibred non-torus 2-bridge knot has 
the universal fibration property.
\vspace{0.1in}

The proof of the Main Theorems utilizes a construction 
initiated in \cite{BW:2}, based on work of \cite{Gol} (and, ultimately, 
the work of Alexander). 
This involves taking a branched cover of $S^3$ over an unknot and lifting 
a branched surface constructed from a train track. 

The outline of the paper is as follows. In \S\ref{sec_back}, we review 
basic definitions and tools needed for the proofs. In \S\ref{sec_family}, 
we introduce a family of pseudo-Anosov braid types and present their 
train tracks. We proceed with the general template constructions in 
\S\ref{sec_templates}, and identify universal subtemplates
to prove the Main Theorem in \S\ref{sec_maintheorems}. 
The final section, \S\ref{sec_conc}, gives a summary of to 
what extent we can classify the universally fibred knots and links.  
We propose a class of examples of pseudo-Anosov
fibrations which we believe do not induce all knot and link types in 
the transverse flow.

%

\section{Background}
\label{sec_back}

Several of the topics below require entire books to do the material
justice: we limit ourselves to the minimal set of ideas necessary
to follow the remaining sections. 

\subsection{Fibred knots and links}

A knot/link $K$ is said to be \df{fibred} if there is a 
fibration $\pi:S^3-K\to S^1$ with fibre a spanning surface 
$\Sigma$ for $K$. The complement of a neighborhood of $K$ in $S^3$ 
is homeomorphic to the mapping torus of a homeomorphism 
$\Phi:\Sigma\to\Sigma$: i.e., 
$S^3-N(K)\cong(\Sigma\times[0,1])/{(x,0)=(\Phi(x),1)}$. 
The spanning surface $\Sigma$ is called the {\bf fibre surface}, and 
the homeomorphism $\Phi$ which is well-defined up to isotopy is called 
the \df{monodromy}. The union of the closed orbits of the 
suspension flow by $\Phi$ on $S^3-N(K)$ 
is called a \df{planetary link} of $\Phi$. 

It follows from the Thurston classification theorem for surface 
homeomorphisms \cite{Thu} that the monodromy of a fibred knot/link 
is isotopic to $\Phi'$, one of the following three types of maps:
\begin{enumerate}
\item {\bf periodic:} $(\Phi')^n=\id$ for some $n$;
\item {\bf pseudo-Anosov:} $\Phi'$ possesses a pair of transverse 
(singular) foliations $\lambda^u$ and $\lambda^s$ along which 
$\Phi'$ is uniformly expanding and contracting respectively;
\item {\bf reducible:} there exists a collection of disjoint simple 
closed curves $\CC \subset\Sigma$ such that 
$\Phi'(\CC) = \CC$, and each connected component of 
$\Sigma-\CC$  has negative Euler characteristic. 
\end{enumerate}
We say that a fibred knot/link $K$ is of periodic, pseudo-Anosov, 
or reducible type respectively if its monodromy map is isotopic to a
periodic, pseudo-Anosov, or reducible map. 


We will be most interested exclusively in those fibred knots/links 
of pseudo-Anosov type (e.g., the figure-eight knot) since 
these alone may permit universal fibrations. 
We may reduce the problem to the case where the monodromy is 
the pseudo-Anosov representative thanks to a theorem 
of Asimov and Franks \cite{AF} which guarantees that any map isotopic to a 
pseudo-Anosov map (and hence any vector field transverse to the 
fibration) may not remove any of these periodic orbits --- isotoping 
the monodromy may generate more, but not less. Hence, the set of 
planetary links of the pseudo-Anosov map is ``minimal'' with regards 
to changing the monodromy map (or the vector field). Henceforth, 
by {\em the} planetary link of a knot/link is meant the planetary
link of its pseudo-Anosov monodromy map. 

\subsection{Braids and branched covers}

We use the language of braids throughout 
the remainder of the paper. Recall that \df{braids} are isotopy classes 
of disjointly embedded arcs monotonically connecting fixed endpoints 
(as in Fig.~\ref{fig_braid}). Braids on $n$ strands form a group
$B_n$ under concatenation, with standard generators $\sigma_i$ 
denoting an elementary crossing of the $\ith$ over the $(i+1)^{st}$ 
strand. Inverses correspond to reversing the 
crossing of the sign.\footnote{There is a tradition of dynamicists
and topologists using opposite sign conventions for braids: we 
employ the dynamicists' convention and apologize for the 
inevitable annoyance.} A \df{closed} or \df{geometric braid} 
is obtained by joining the two sets of endpoints around some fixed
\df{braid axis}: algebraically, this corresponds to taking the 
conjugacy class of the braid in the braid group. 

The classification of monodromies into periodic, pseudo-Anosov, or
reducible type has its analogue for braids. 
Let $\beta$ denote the closure of $b$, an $n$-braid with braid axis 
$\alpha$. Since $\alpha$ is unknotted, its complement fibres over
$S^1$ with fibre a disc $A$ meeting $\beta$ at $n$ points. 
The monodromy map on $A$ is fixed on the boundary and, by 
removing the intersection with $\beta$, becomes a map on an
$n$-punctured disc. The Thurston Classification theorem applied
to this monodromy map implies that one may classify a braid $b$
as being of {\bf periodic type}, {\bf pseudo-Anosov type}, or 
{\bf reducible type}.  
The pseudo-Anosov type braids will be of central importance in 
analyzing fibrations. 

While fibrations of knot complements are difficult to visualize
(the explanations of the trefoil, in \cite[10I]{Rolfsen}, and the 
figure-eight knot, in \cite{BW:2} are not short!), the \df{branched 
covering} construction of Goldsmith \cite{Gol} and Birman 
\cite{Bir} provides a concrete mechanism for analysis. We use the 
method as presented in \cite[pp. 26-30]{BW:2}, keeping similar notation. 

\begin{lem}\cite[Lemma 1]{Gol}
\label{lem_goldsmith}
Let $\beta$ be a single-component unknotted geometric braid with 
braid axis $\alpha$ and meridional disc $A$ spanning $\alpha$. 
Let $p:S^3\rightarrow S^3$ be a branched covering 
space projection whose branch set is $\beta$.  
Then $\underline{\alpha} = p\inv(\alpha)$ is a non-trivial 
fibred knot/link in $S^3$ with fibre $\underline{A} = p\inv(A)$. 
\end{lem}

A result of Birman \cite{Bir} shows that every fibred knot may be obtained
as $p\inv(\alpha)$ for some $\beta$ and some (irregular) cover $p$. 
Moreover, the fibration of the complement of $\underline{\alpha}$ 
is itself a lift of the fibration of the complement of $\alpha$.

For the remainder of this paper, we work with pairs of curves
$(\alpha,\beta)$ which are \df{exchangeable},
i.e., there is an isotopy of $S^3$ which presents
$\alpha$ as a geometric $n$-braid with braid axis $\beta$ and vice versa
\cite{MR}. Let $a$ (resp. $b$) denote the braid in $B_n$ whose closure
is $\alpha$, with $\beta$ as its axis (resp. the braid whose closure
is $\beta$ with axis $\alpha$). In this case, $\underline{\alpha}$ 
is the closure of $a^m$ if $p$ is $m$-fold.

\subsection{Train tracks and templates}

The problem of how to understand the complex dynamics associated to 
a pseudo-Anosov surface homeomorphism was resolved with the theory 
of \df{train tracks}. The central idea, going back to Williams 
\cite{W}, is to project the foliation down to a branched manifold. 
We call attention to \cite{PH,BH} for treatments of this 
theory from different perspectives. 

Recall that a pseudo-Anosov surface map is characterized by the pair of 
transverse singular foliations, $\lambda^s$ and $\lambda^u$, along which 
the map is uniformly contracting and expanding respectively. 
Roughly speaking, a \df{train track} for a pseduo-Anosov map 
is a branched 1-manifold obtained by cutting open the singular
foliation $\lambda^s$ in a certain manner and collapsing each leaf 
to a point: see \cite{GK} for examples relevant to this paper.
The dynamics of the pseduo-Anosov homeomorphism is then easily represented 
by an induced self-immersion of the train track. We note that there are 
several presentations available for train tracks. For the 
remainder of this paper, we follow e.g., \cite{BH,CH} and represent 
the train track as a graph. The induced map is then easily
represented as a map taking each edge of the graph to a sequence
of edges. Two simple examples of train track graphs and associated
graph maps are illustrated in Fig.~\ref{fig_graph}.


The problem of how to capture all of the knot and link data of a 
dynamically complex 3-dimensional flow was successfully tackled 
by Birman and Williams in the early 1980s \cite{BW:1,BW:2} in the 
theory of templates (a.k.a. ``knotholders''). 

\begin{figure}[hbt]
\begin{center}
\psfragscanon
\psfrag{1}[][]{\Large $\x_1$}
\psfrag{2}[][]{\Large $\x_2$}
\psfrag{3}[][]{\Large $\x_3$}
\psfrag{4}[][]{\Large $\x_4$}
\includegraphics[angle=0,width=4.5in]{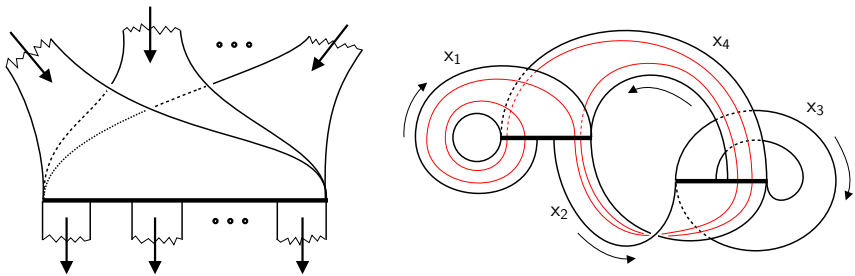}
\caption{A branchline chart for a template [left]; 
an example of a template with an unknotted periodic orbit on it [right].
In the alphabet $\{\x_1,\x_2,\x_3,\x_4\}$, this knot corresponds to the 
word $\word{\x_1^2(\x_2\x_4)^2}=\word{\x_1\x_1\x_2\x_4\x_2\x_4}$.}
\label{fig_charts}
\end{center}
\end{figure}

\begin{dfn}
\label{def_Template}
A \df{template} is an embedded branched 2-manifold with boundary,
outfitted with an expansive semiflow. Templates have a natural 
decomposition into a finite number of \df{branchline charts} 
(Fig.~\ref{fig_charts}[left]), each containing a number ($>1$) of 
incoming and outgoing strips whose ends are joined respecting the 
direction of the semiflow.
\end{dfn}

The term ``semiflow'' means that the flow is well-defined in forward time,
but not in backward time (at the branch lines in particular --- there
is a loss of uniqueness). By ``expansive,'' it is meant that the semiflow
expands volume everywhere. 

Every template carries an infinite link of periodic orbits of the 
semiflow. There is a well-defined symbolic language for describing these
orbits. Upon labeling each strip of a template $\T$ with a symbol
$\x_i$, it is clear that to each forward orbit of the semiflow on $\T$ is
associated a semi-infinite word in the alphabet $\{\x_i\}$. Furthermore,
from the expansivity of the semiflow, it can be shown that periodic 
flowlines on $\T$ are in fact in bijective correspondence 
with (admissible) periodic words in this alphabet, 
up to cyclic permutations of the words (see \cite{BW:2} or 
\cite[Lem. 2.4.1]{GHS} for details). 

The Template Theorem \cite[Theorem 2.1]{BW:2}
of Birman and Williams gives very general conditions under which 
the dynamics of a three-dimensional flow is accurately captured by a
template: see \cite{GHS} for a comprehensive treatment. 
For our purposes, we note that for a fibred knot/link $K$ with 
pseudo-Anosov monodromy, there is a template $\T \subset S^3-K$ 
whose knot and link types are in bijective correspondence with those 
of the flow except for a finite number of orbits on the boundary
of the template. The template is obtained by suspending the 
train track graph map for the monodromy, and cutting along the 
vertices of the graph as necessary. These cuts may 
produce ``extra'' orbits on the templates not in the original flow. To 
minimize confusion, we note that for the remainder of this 
paper, any such exceptional orbits which we ``lose'' are unknots,
and are easily found elsewhere in the interior of the template. 

\begin{dfn}
A \df{universal template} in $S^3$ is one which possesses all knot and 
link types as closed orbits of the semiflow. 
\end{dfn}

The paper \cite{G} established that universal templates exist
in abundance. From this work and \cite{GHS} comes a 
practical criterion for determining when a given template
is universal. One additional definition is needed. Given any
closed orbit $\kappa$ on a template $\T\subset S^3$, the 
\df{twist}, $\tau(\kappa)$, is defined to be the twist
number of the normal bundle to $\T$ along $\kappa$ (which is 
either an annulus or M\"obius band having $\kappa$ as its core). 
For the remainder of this paper, we consider twists for 
unknotted curves only, avoiding any ambiguity in how to count
twist: an unknot with zero twist is one whose normal bundle 
can be isotoped to be a plane annulus.

The following criterion is a very slight modification of 
\cite[Cor. 3.2.17]{GHS}, with a near-identical proof. Since
the method of proof relies heavily on terminology and 
techniques of \cite{GHS}, we suppress the details. 

\begin{figure}[hbt]
\begin{center}
\psfragscanon
\psfrag{a}[][]{{\Huge $\kappa$}}
\psfrag{b}[][]{{\Huge $\kappa'$}}
\psfrag{c}[][]{{\Huge $\kappa''$}}
\psfrag{x}[][]{{\Large $+$}}
\psfrag{z}[][]{{\Large $-$}}
\includegraphics[angle=0,width=5.0in]{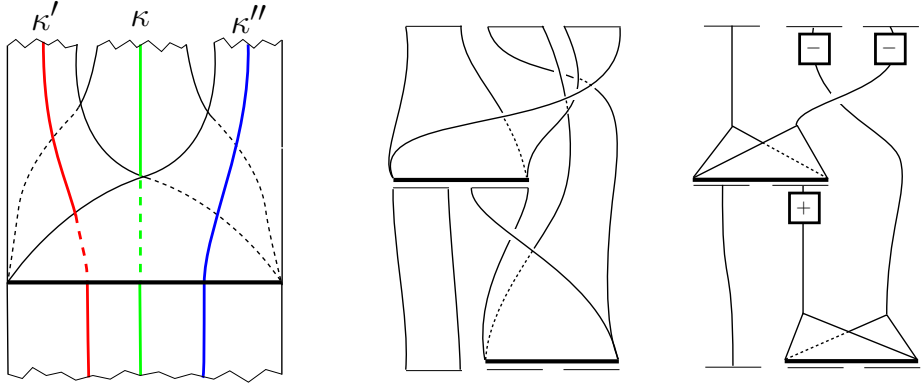}
\caption{$\kappa', \kappa, \kappa''$ from left to right [left];
a template and its cartoon representation --- top and bottom 
are identified [center,right].}
\label{fig_UT}
\label{fig_three_orbits}
\end{center}
\end{figure}

\begin{thm}
\label{thm_UT}
Let ${\mathcal T}$ be a template in $S^3$. 
Suppose that there exist three disjoint closed orbits on ${\mathcal T}$, 
$\kappa$, $\kappa'$, $\kappa''$ such that
\begin{enumerate}
\item they are separable unlinked unknots; 
\item $\tau(\kappa) = 0$, $\tau(\kappa')>0$, $\tau(\kappa'')<0$; and 
\item these three unknots intersect some branchline of $\T$ as in 
Fig.~\ref{fig_three_orbits}[left] with the specified adjacencies and 
strip crossings. 
\end{enumerate}
Then $\T$ is a universal template.
\end{thm}

The proofs of all the main theorems in this paper consist of 
constructing and simplifying templates for fibred knots and links, and then
searching for the three unknots $\{\kappa, \kappa', \kappa''\}$ to which 
one may apply Theorem~\ref{thm_UT}. Because these templates are
quite complex, we introduce a ``cartoon'' notation (cf. \cite{Tuff}).
\begin{enumerate}
\item Strips in the template are collapsed along the direction transverse
  to the semiflow to become braided curves, except in a neighborhood of the 
  branch lines;
\item Each half-twist of each strip is represented by a 
  label $+$ or $-$ within a box, depending on the sign of the crossing;
\item Templates will often be ``cut open'' with the top and bottom
  to be identified.
\end{enumerate}
For example, the template in Fig.~\ref{fig_UT}[center]  
is expressed in cartoon form in Fig.~\ref{fig_UT}[right]. 

\subsection{Templates for planetary links}

The branched cover construction allows one to construct a template
for the planetary link of certain fibred knots/links. Recall the
notation that $\beta$, the closure of $b$, is an unknotted geometric 
braid with braid axis $\alpha$; $p$ is an $m$-fold branched cover over 
$\beta$; and $\underline{\alpha}=p\inv(\alpha)$ is the lift of the 
braid axis, which, in the case of exchangeable braids is the closure
of the braid $a^m$. 

We summarize the procedure of \cite[pp. 26-30]{BW:2} in the 
case of an exchangeable braid. If $b$ is a braid of pseudo-Anosov type, 
then $\underline{\alpha}$ is fibred and also of pseudo-Anosov type, 
since the monodromy map for $\underline{\alpha}$ projects under $p$ to the 
monodromy map of $\alpha$ represented by the pseudo-Anosov braid type $b$. 

\begin{figure}[hbt]
\begin{center}
\psfragscanon
\psfrag{+}[][]{$+$}
\psfrag{-}[][]{$-$}
\includegraphics[angle=0,width=5.0in]{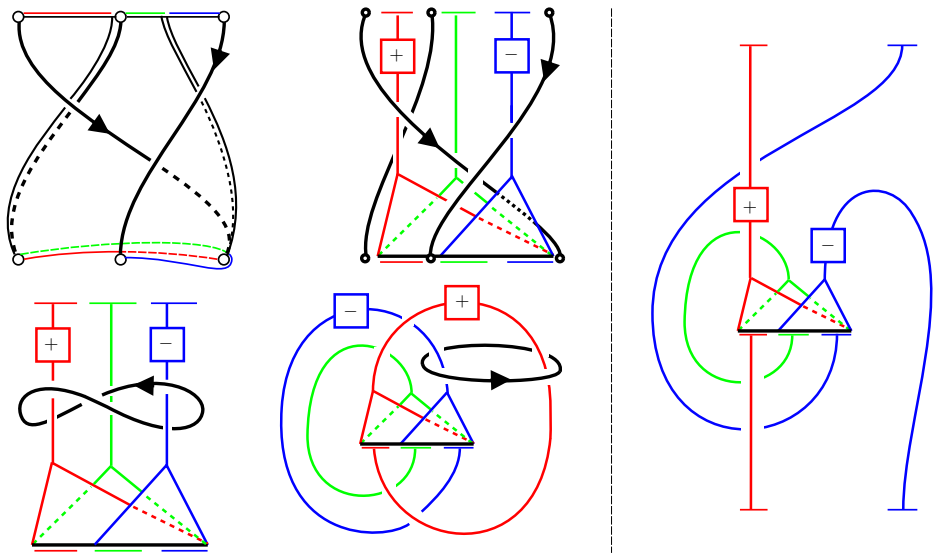}
\caption{The branched cover construction for the figure-eight
knot complement. The suspension of the train track for the 
braid $b=\sigma_1\sigma_2^{-1}$ is a template [upper left]. By peeling off
a copy of $\beta$ and cutting along its spanning disc [lower left],  
one has a fundamental domain for the template of the $m$-fold cover [right].}
\label{fig_8ex}
\end{center}
\end{figure}

The following procedure yields a template for the fibration 
of the closure of $a^m$: 
\begin{enumerate}
\item Construct a train track graph for $b$ along with 
its induced map; 
\item Suspend the train track graph map 
and cut as necessary to obtain a template $\T$ in the solid torus 
$S^3-\alpha$; 
\item Peel off a copy of the unknot $\beta$ from the boundary 
of $\T$ and find a spanning disc $D$ for this unknot; 
\item Cut $\T$ along 
$D$ and glue $m$ copies of this cut template together end-to-end cyclically 
to obtain $\U^m$, the template for the fibration of the 
complement of $\underline\alpha$. 
\end{enumerate}
This procedure was carried out carefully in \cite[pp. 28-30]{BW:2}
for the figure-eight knot, which is the closure of the braid
$(\sigma_1\sigma_2\inv)^2$. Note that in this case 
$a=b=\sigma_1\sigma_2\inv$ is the simplest pseudo-Anosov braid 
type. The train track graph for this braid is an interval, and the 
graph map is represented in Fig.~\ref{fig_graph}[left]. 
We reproduce this example in cartoon notation in Fig.~\ref{fig_8ex} 
in preparation for the proof of Main Theorem (c). Familiarity with
this example will reveal patterns in the general case. 

%

\section{Families of fibred knots and links} 
\label{sec_family}

We focus on the following families of braids along which 
to perform the branched cover construction. 

\begin{dfn}
Given integers $n_1, n_2, \cdots, n_k \ge 1$ (with $k>1$), denote by 
$b:=\braid$ the braid on $(1+\sum_{r = 1}^{k}n_r)$ strands given by 
the braid word 
\begin{equation}
\label{eq_b}
\underbrace{
\sigma_1 \sigma_2 \cdots \sigma_{n_1}}_{n_1} 
\underbrace{
\sigma_{n_1 + 1}^{-1} \cdots \sigma_{n_1 + n_2}^{-1}}_{n_2} 
\cdots 
\underbrace{
\sigma_{1+\sum_{r = 1}^{k-1}n_r}^{{(-1)}^{(k+1)}}
 \cdots 
\sigma_{\sum_{r = 1}^{k}n_r}^{{(-1)}^{(k+1)}}}_{n_k}.
\end{equation}
\end{dfn}

Thus, $b$ is equivalent to the form given in (\ref{eq_Main}). 
For example, the braid in Fig.~\ref{fig_braid} is $b_{(2,2,1,1)}$. 
Clearly, one can also use the mirror image of these braids, and all the 
theorems we prove hold here as well.


\begin{lem}
\label{lem_exchange}
The closure of $\braid$ is self-exchangeable. Specifically, given 
$\beta$ the closure of $\braid$ with braid axis $\alpha$, there exists
an isotopy of $S^3$ taking $\alpha$ to the closure of $\braid$ 
with braid axis $\beta$.  
\end{lem}

\begin{figure}[hbt]
\begin{center}
\psfragscanon
\psfrag{b}[][]{\Huge $b$}
\includegraphics[angle=0,width=5in]{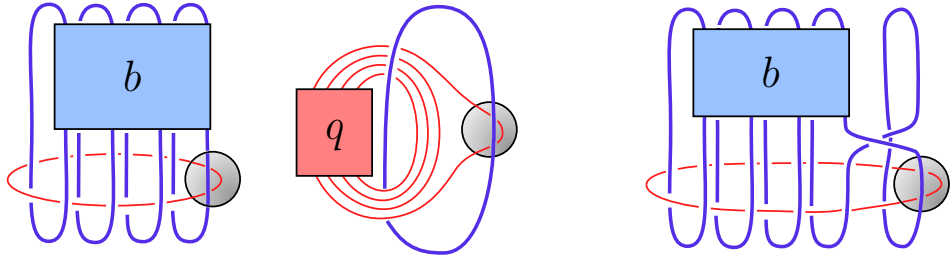}
\caption{The induction hypothesis: the closure of the braid $b$ can be
exchanged with its axis rel a fixed ball [left]; the setup for 
the induction step [right].}
\label{fig_Exchange1}
\end{center}
\end{figure}

\begin{figure}[hbt]
\begin{center}
\psfragscanon
\psfrag{1}[][]{}
\psfrag{2}[c][c]{\rm twist}
\psfrag{3}[c][c]{\rm flip}
\psfrag{4}[c][c]{\rm slide}
\psfrag{b}[c][c]{\Huge $b$}
\includegraphics[angle=0,width=5.0in]{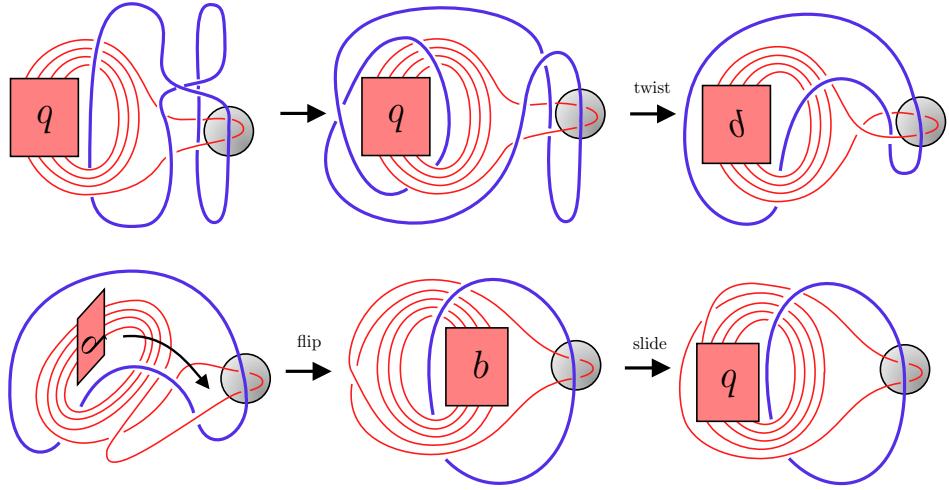}
\caption{The induction step: the closure of $b\sigma_N$ can be
exchanged with its axis rel a fixed ball.}
\label{fig_Exchange2}
\end{center}
\end{figure}

\pf
We induct on the length of the braid word $b$. Assume as an induction 
hypothesis that the $N$-strand closed braid $\beta$ exchanges with 
$\alpha$ as in Fig.~\ref{fig_Exchange1}[left]: specifically, that there 
is an isotopy which is fixed on a solid ball as indicated. Now,
for a braid of the form $b\sigma_{N}$ (in Fig.~\ref{fig_Exchange1}[right])
or $b\sigma_{N}\inv$, fit the final strand within a larger fixed solid ball 
and perform the isotopy given by induction. Fig.~\ref{fig_Exchange2} shows
that this final strand of $\beta$ may be exchanged relative to 
a yet smaller fixed ball. Note that in the diagram, the set of 
crossings in the braid word $b$ is flipped and rotated, but is kept 
together rigidly.
\qed

\begin{lem}
\label{lem_pA}
The braid $b=\braid$ is of pseudo-Anosov type.
\end{lem}
A proof of this can be explicitly given by one of the [several] 
algorithms available.\footnote{We have used the Bestvina-Handel algorithm 
\cite{BH} for this family: though straightforward, the algorithm is 
lengthy and of little consequence for the remainder of this article.} 
We note, however, that the braid $b$ is precisely that 
considered by Gabai and Kazez \cite{GK}, who show the pseudo-Anosov 
braid type. 

By Lemma~\ref{lem_pA}, we have the following; 
\begin{cor}
For each integer $m > 1$, the closure of $(\braid)^m$ 
is a fibred knot/link of pseudo-Anosov type. 
\end{cor}

\begin{figure}[htbp]
\begin{center}
\psfragscanon
\psfrag{a}[][]{\Large $e_i^0$}
\psfrag{d}[][]{\Large $e_i^0$}
\psfrag{e}[][]{\Large $e_i^1$}
\psfrag{g}[][]{\Large $e_i^{n_i-1}$}
\psfrag{h}[][]{\Large $e_i^{n_i}$}
\includegraphics[angle=0,width=4in]{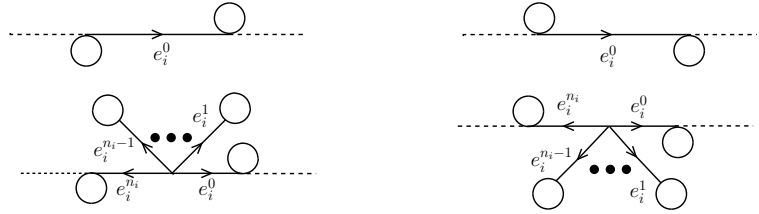}
\caption{
$i$ odd, $n_i=1$ [upper left]; 
$i$ even, $n_i=1$ [upper right];
$i$ odd, $n_i >1$ [lower left];
$i$ even, $n_i >1$ [lower right].}
\label{fig_graph_def}
\end{center}
\end{figure}

\begin{figure}[hbt]
\begin{center}
\psfragscanon
\psfrag{a}[][]{\Large $e_1^0$}
\psfrag{b}[][]{\Large $e_2^2$}
\psfrag{c}[][]{\Large $e_2^0$}
\psfrag{d}[][]{\Large $e_2^1$}
\psfrag{e}[][]{\Large $e_1^0$}
\psfrag{f}[][]{\Large $e_2^0$}
\psfrag{g}[][]{\Large $e_3^0$}
\psfrag{h}[][]{\Large $e_1^2$}
\psfrag{i}[][]{\Large $e_1^1$}
\psfrag{j}[][]{\Large $e_1^0$}
\psfrag{k}[][]{\Large $e_2^0$}
\psfrag{l}[][]{\Large $e_3^3$}
\psfrag{m}[][]{\Large $e_3^2$}
\psfrag{n}[][]{\Large $e_3^1$}
\psfrag{o}[][]{\Large $e_3^0$}
\includegraphics[angle=0,width=4.5in]{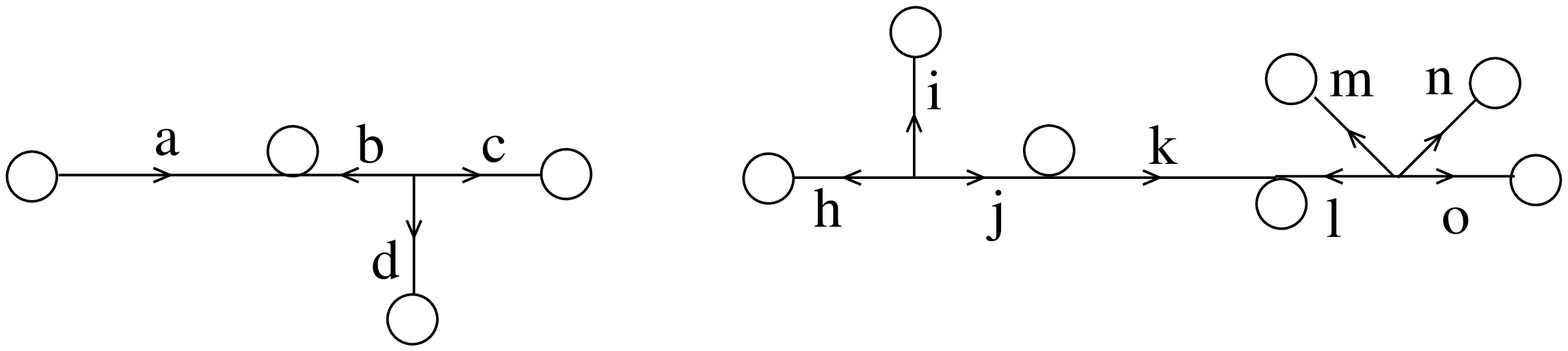}
\caption{$G_{(1,2)}$ [left]; $G_{(2,1,3)}$ [right].}
\label{fig_ex_graph}
\end{center}
\end{figure}

%
%
The braid $\braid$ induces a very simple map on the closed disc 
consisting of a collection of $k$ rotations which alternate direction,
and the associated train track is very straightforward to compute.
Let $G = G_{(n_1, \cdots, n_k)}$ be a chain of $k$ radial ``stars'' 
$S_1, \cdots, S_k$, each star having valence $(n_i+1)$ vertex 
if $n_i > 1$ (Fig.~\ref{fig_graph_def}). 
The radial stars alternate their up-down orientation, and small loops 
(corresponding to the periodic orbit) 
are attached at the ends of the radial oriented edges 
as in Fig.~\ref{fig_graph_def}. 
The set of non-loop edges in $S_i$ equals $\{e_i^0\}$ if $n_i=1$ 
(resp. $\{e_i^0, \cdots, e_i^{n_i}\}$ if $n_i > 1$). 
For example, see Fig.~\ref{fig_ex_graph}.

\begin{figure}[hbt]
\begin{center}
\psfragscanon
\psfrag{X}[][]{\huge $\Rightarrow$}
\psfrag{b}[][]{\huge $e_i^0$}
\psfrag{B}[][]{\huge $g_i(e_i^0)$}
\psfrag{p}[][]{\huge $e_i^0$}
\psfrag{r}[][]{\huge $e_i^{n_i-1}$}
\psfrag{s}[][]{\huge $e_i^{n_i}$}
\psfrag{P}[l][l]{\huge $g_i(e_i^0)$}
\psfrag{R}[r][r]{\huge $g_i(e_i^{n_i-1})$}
\psfrag{S}[l][l]{\huge $g_i(e_i^{n_i})$}
\includegraphics[angle=0,width=5.0in]{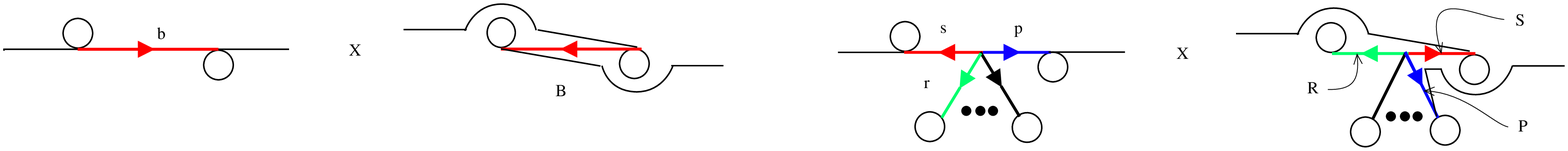}
\caption{$i$ even, $n_i = 1$ [left]; 
$i$ even, $n_i > 1$ [right].}
\label{fig_graphmap}
\end{center}
\end{figure}

\begin{lem}
\label{lem_traintrack}
An efficient 
\footnote{For the definition of the {\bf efficient graph maps}, 
see \cite[pp. 114]{BH}. In general, the efficient graph map is not 
unique for the braid $b$.} graph map
$g = g_{(n_1, \cdots, n_k)}:G \rightarrow G$ for 
the train track graph of $b_{(n_1, \cdots, n_k)}$ 
is given as $g_k \circ \cdots \circ g_1$, 
where $g_i$ rotates the $i$-th star $S_i$ in the 
clockwise (resp. counterclockwise) direction if 
$i$ is even (resp. odd), see Fig.~\ref{fig_graphmap}. 
\end{lem}


The algebraic presentation of the train track map is simple to 
derive knowing that it is composed of rotations. For 
simplicity, we ignore all loop-edges in the train track. 
Given an oriented edge $E$, we denote the same edge oriented in the 
opposite direction by $\overline{E}$. 
The following lemma is used to derive a normal form of the template 
$\T_{(p,q)}$ in Section~\ref{subsec_NormalForm}.

\begin{figure}[hbt]
\begin{center}
\psfragscanon
\psfrag{x}[][]{$\Rightarrow$}
\includegraphics[angle=0,width=5in]{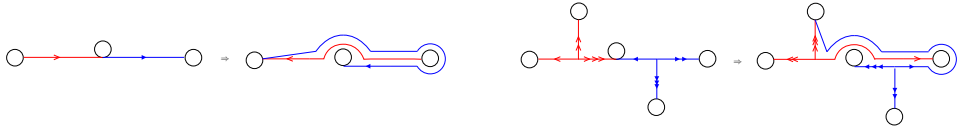}
\caption{$g_{(1,1)}$ [left]; $g_{(2,2)}$ [right].}
\label{fig_graph}
\end{center}
\end{figure}

\begin{lem}
\label{lem_graph}
The images of the edges under $g=g_{(p, q)}$ are as follows; 
For $p = q = 1$, 
\begin{eqnarray*}
g(e_1)&=& \overline{e_2^0}~\overline{e_1^0}, 
\\ 
g(e_2)&=& e_1^0 e_2^0 \overline{e_2^0}. 
\end{eqnarray*}
For $p, q >1$, 
\begin{eqnarray*}
g(e_1^{p}) &=& e_1^0 \overline{e_2^{q}} e_2^0,
\\
g(e_2^{q}) &=&e_2^0 \overline{e_2^0} e_2^{q} 
\overline{e_1^0}e_1^1,
\\
g(e_1^i) &=& e_1^{i+1}\ 
(i = 0, \cdots, p -1),
\\
g(e_2^i) &=& e_2^{i+1}\ 
(i = 0, \cdots, q -1), 
\end{eqnarray*}
see Fig.~\ref{fig_graph}.
\end{lem}

%

\section{Subtemplates in branched covers}
\label{sec_templates}

The various cases of the Main Theorem all possess the same 
strategy of proof. Given the braid $b=\braid$, we suspend the 
train track $G=G_{(n_1,\ldots,n_k)}$ 
by the train track map $g=g_{(n_1,\ldots,n_k)}$
to obtain a template $\T_{(n_1,\ldots,n_k)}$.
We then use the branched covering method to obtain a template 
$\U_{(n_1,\ldots,n_k)}^m$ which captures the monodromy of the 
closure of the braid $b^m$. Since, in general, $\U_{(n_1,\ldots,n_k)}^m$ 
is much too complex to visualize, we focus on a particular 
subtemplate which we may show to be universal. In this section, 
we prove the technical results for restricting attention to 
appropriately simple subtemplates. 

\subsection{A normal form}
\label{subsec_NormalForm}

We now derive the branched surface $\T_{(p,q)}$ obtained 
from $G_{(p,q)}$ by suspending $g_{(p,q)}$  from
the previous section. In what follows, we assume $p,q>1$. 
(Other cases are much simpler and left to readers.) 
Let $q_{\l}$ and $q_r$ are the valence $(p+1)$ vertex and $(q+1)$ vertex 
of $S_1$ and $S_2$ respectively. 
These suspend to give the knots $k_{\l}$ and $k_r$. 
Notice that the set of loop-edges in $S_1 \cup S_2$ 
is a periodic orbit, say $A_{(p,q)}$, 
corresponding to the closure of $b_{(p,q)}$.

\begin{figure}[hbt]
\begin{center}
\psfragscanon
\psfrag{a}[][]{}
\psfrag{b}[][]{}
\psfrag{e}[][]{$e_1^p$}
\psfrag{f}[][]{$e_1^{p-1}$}
\psfrag{g}[][]{$e_1^2$}
\psfrag{h}[][]{$e_1^1$}
\psfrag{i}[][]{$e_1^0$}
\psfrag{j}[][]{$e_2^q$}
\psfrag{k}[][]{$e_2^{q-1}$}
\psfrag{l}[][]{$e_2^2$}
\psfrag{m}[][]{$e_2^1$}
\psfrag{n}[][]{$e_2^0$}
\psfrag{x}[][]{\Large $\Rightarrow$}
\includegraphics[angle=0,width=5in]{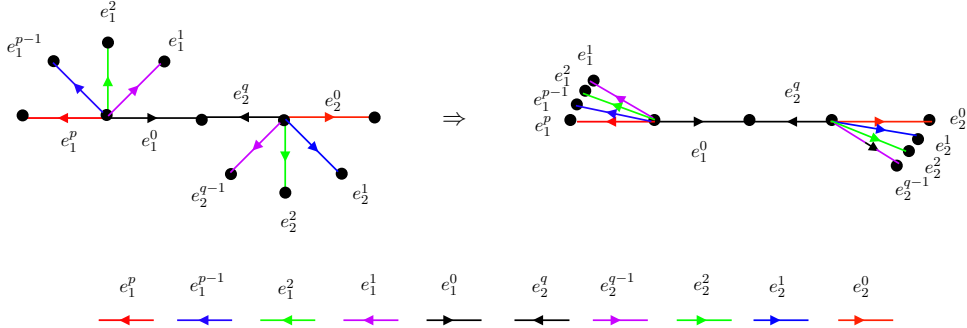}
\caption{Flattening a pair of stars into strips.}
\label{fig_split}
\end{center}
\end{figure}

To visualize $\T_{(p,q)}$ as a template, we need to ``cut open'' 
$q_{\l}$, $q_r$ and $A_{(p,q)}$, 
and flatten out the strips in some canonical
way. This is illustrated in Fig.~\ref{fig_split}. 

Recall that $g(e_1^p)=e_1^0 \overline{e_2^q}e_2^0$.
We split the edge $e_1^p$ into $E$ and $E'$ so that
$g(E) = e_1^0 \overline{e_2^q}$ and $g(E') = e_2^0$.
Split the edge $e_2^q$ into $F, F'$ and $F''$ so that
$g(F) = e_2^0$, $g(F')= \overline{e_2^0}$ and
$g(F'') = e_2^{q}\overline{e_1^0}e_1^1$.
For the template $\T_{(p,q)}$, we can use the symbols $\x_{i=0..p+1}$,
and $\y$,$\z_{j=0...q+1}$ for the first star, and second 
star respectively. The edge $e_1^0$ goes to $\x_0$, and  
the edge $e_2^0$ goes to $\z_0$, etc., and $e_1^p$ is split
into $\x_p$ and $\x_{p+1}$ corresponding to $E$ and $E'$, and
$e_2^q$ split into $\z_q,\z_{q+1},\y$ corresponding to $F,F'$ and
$F''$.
We denote the left end point and right end point
of the branch segment of the strip 
$\x$ by $p_{\l}(\x)$ and $p_r(\x)$ respectively. 
Since the initial vertex of 
$e_1^0, e_1^1, \cdots, e_1^p$ is the same vertex $q_{\l}$, 
we must abstractly identify 
(1) $p_r(\x_1)$, $p_r(\x_2)$ ,$\cdots$, $p_r(\x_{p})$ and 
$p_{\l}(\x_0)$. In the same manner, we need to identify 
(2) $p_r(\x_0)$ and $p_{\l}(\y)$, 
(3) $p_r(\y)$ and $p_{\l}(\z_{q+1})$, 
(4) $p_r(\z_{q+1})$ and $p_{\l}(\z_q)$, 
(5) $p_{\l}(\z_0)$, $p_{\l}(\z_1)$, $\cdots$, $p_{\l}(\z_{q-1})$ 
and $p_r(\z_q)$, and 
(6) $p_{\l}(\x_p)$ and $p_{r}(\x_{p+1})$. 

For the remainder of the paper, we conveniently leave these
boundary curves as distinct. When constructing subtemplates of
this branched surface, we are careful never to use these
edges within the subtemplates. In this way, we do not need
to worry about the changes to the periodic orbits set we have
made. 

\begin{figure}[htbp]
\begin{center}
\psfragscanon
\psfrag{1}[][]{$\x_{p+1}$}
\psfrag{2}[][]{$\x_p$}
\psfrag{3}[][]{$\x_{p-1}$}
\psfrag{4}[][]{$\x_{3}$}
\psfrag{5}[][]{$\x_1$}
\psfrag{6}[][]{$\x_0$}
\psfrag{7}[][]{$\y$}
\psfrag{8}[][]{$\z_{q+1}$}
\psfrag{9}[][]{$\z_{q}$}
\psfrag{10}[][]{$\z_{q-1}$}
\psfrag{11}[][]{$\z_3$}
\psfrag{12}[][]{$\z_1$}
\psfrag{13}[][]{$\z_0$}
\psfrag{14}[][]{$\x_2$}
\psfrag{15}[][]{$\z_2$}
\psfrag{a}[][]{$+$}
\psfrag{c}[][]{$-$}
\includegraphics[angle=0,width=4in]{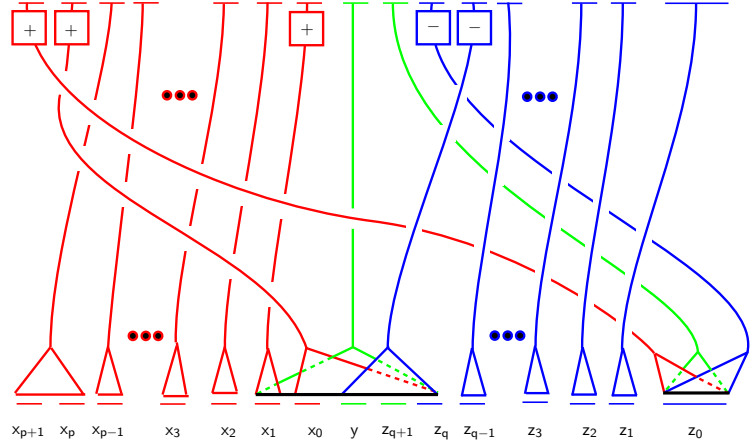}
\caption{The normal form for the template $\T_{(p,q)}$.}
\label{fig_normal_form2}
\end{center}
\end{figure}

The normal form of $\T_{(p,q)}$ appears in cartoon form 
as in Fig.~\ref{fig_normal_form2}. 

In order to construct $\U^m_{(p,q)}$, the template for the $m$-fold 
branched cover of $\T_{(p,q)}$, 
we must ``peel off'' a parallel copy of $\beta$
and cut $\T_{(p,q)}$ by a spanning disc for this unknot: cf. 
the figure-eight knot case in Fig.~\ref{fig_8ex}. Thus, 
the fundamental domain of $\U^m_{(p,q)}$ depends on the linking 
properties of $\beta$ with $\T_{(p,q)}$. 

\begin{figure}[htbp]
\begin{center}
\psfragscanon
\psfrag{1}[][]{$\x_{p+1}$}
\psfrag{2}[][]{$\x_p$}
\psfrag{3}[][]{$\x_{p-1}$}
\psfrag{4}[][]{$\x_{3}$}
\psfrag{5}[][]{$\x_1$}
\psfrag{6}[][]{$\x_0$}
\psfrag{7}[][]{$\y$}
\psfrag{8}[][]{$\z_{q+1}$}
\psfrag{9}[][]{$\z_{q}$}
\psfrag{10}[][]{$\z_{q-1}$}
\psfrag{11}[][]{$\z_3$}
\psfrag{12}[][]{$\z_1$}
\psfrag{13}[][]{$\z_0$}
\psfrag{14}[][]{$\x_2$}
\psfrag{15}[][]{$\z_2$}
\psfrag{a}[][]{$+$}
\psfrag{c}[][]{$-$}
\includegraphics[angle=0,width=4in]{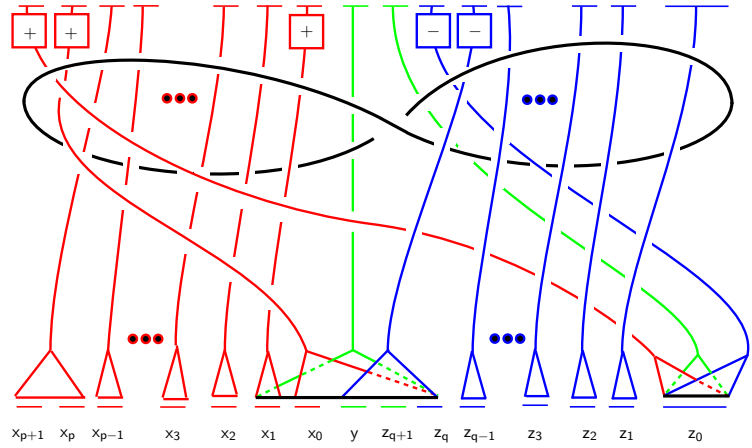}
\caption{A parallel copy of the closed braid $\beta$ links 
certain strips of $\T_{(p,q)}$.}
\label{fig_normal_form3}
\end{center}
\end{figure}

\begin{figure}[htbp]
\begin{center}
\psfragscanon
\psfrag{a}[][]{}
\psfrag{b}[][]{}
\psfrag{x}[][]{\huge $k_{\l}$}
\psfrag{y}[][]{\huge $k_c$}
\psfrag{z}[][]{\huge $k_r$}
\psfrag{w}[][]{\huge $\beta$}
\psfrag{X}[][]{\huge $k_{\l}'$}
\psfrag{Y}[][]{\huge $k_c'$}
\psfrag{Z}[][]{\huge $k_r'$}
\includegraphics[angle=0,width=3.5in]{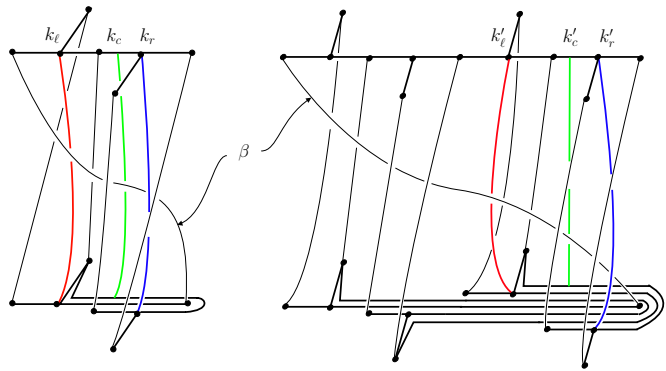}
\caption{$p=q=2$ [left]; $k=4$ and $n_1=n_2=n_3=n_4=2$ [right].}
\label{fig_linking_num4}
\end{center}
\end{figure}

\begin{lem}
\label{lem_linking}
Let $\beta$ be the closure of $b_{(p,q)}$ and $\beta'$ a parallel
copy. Then $\beta'$ links $\T_{(p,q)}$ as in 
Fig.~\ref{fig_normal_form3}. 
\end{lem}
\pf 
Note that $g_{(n_1,n_2)}$ 
has three fixed points $q_{\l}$, $q_c$, $q_r$, 
where $q_{\l}$ and $q_r$ are as above and $q_c \in e_{2}^q$ 
(see Lemma~\ref{lem_graph}). 
These suspend to give the three unknots $k_{\l}$, $k_c$, and $k_r$. 
It is clear from the image of the train track (see
Figs.~\ref{fig_graph} and \ref{fig_linking_num4}) that the fixed 
point $q_c$ is ``in back'' of the image of all the other edges, 
and that consequently the suspension of this fixed point,
$k_{c}$, is an unknot which is of twist zero and separable
from all other orbits on the template. Hence, $\lk(k_{c}, \beta) = 0$.
Since the train track map rotates each star about its central vertex
with opposite orientations, it 
is easy to see that $\lk(k_{\l}, \beta) = +1$ and 
$\lk(k_{r}, \beta) = -1$ (see Fig.~\ref{fig_linking_num4}[left]). 

Since the strip $\y$ contains $k_c$ and $\lk(k_{c}, \beta) = 0$, 
$\y$ and $\beta$ are linked as in Fig.~\ref{fig_linking_num4}. 
Because of the way we unfolded the star in Fig.~\ref{fig_split}, 
the way that $\z_{q+1}$ and $\beta$ are linked is same as that of 
$\y$ and $\beta$. 

Notice that the edges $e_1^0, \cdots, e_1^p$ were all ``split off'' of 
$q_{\l}$, and $\lk(k_{\l}, \beta) = +1$. 
Hence $\x_0, \cdots, \x_p, \x_{p+1}$ and $\beta$ are linked 
as in Fig.~\ref{fig_linking_num4}. 
Since the edges 
$e_2^0, \cdots, e_2^{q-1}, F (\subset e_2^{q}) $ were also ``split off'' of 
$q_{r}$ and $\lk(k_{r}, \beta) = -1$, the strips  
$\z_0, \cdots, \z_{q-1}, \z_q$ and $\beta$ are linked 
as in Fig.~\ref{fig_linking_num4}. This completes the proof. 
\qed


\begin{figure}[htbp]
\begin{center}
\psfragscanon
\psfrag{1}[][]{}
\psfrag{2}[][]{}
\psfrag{3}[][]{}
\psfrag{4}[][]{}
\psfrag{5}[][]{}
\psfrag{6}[][]{}
\psfrag{7}[][]{}
\psfrag{8}[][]{}
\psfrag{9}[][]{}
\psfrag{10}[][]{}
\psfrag{11}[][]{}
\psfrag{12}[][]{}
\psfrag{13}[][]{}
\psfrag{14}[][]{}
\psfrag{15}[][]{}
\psfrag{a}[][]{$+$}
\psfrag{b}[][]{}
\psfrag{c}[][]{$-$}
\includegraphics[angle=0,width=5in]{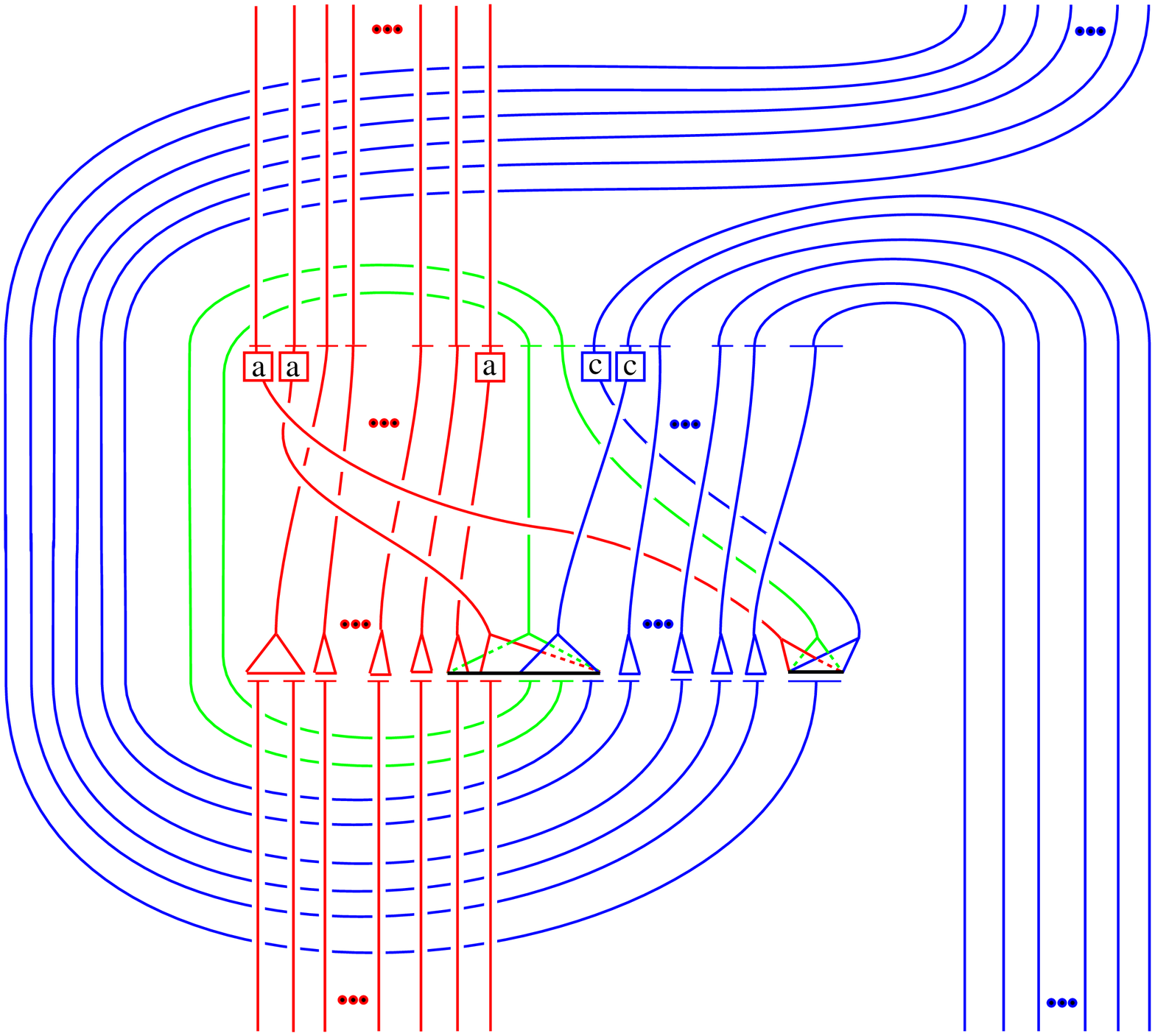}
\caption{The fundamental domain of $\U^m_{(p,q)}$, 
cf. Fig.~\ref{fig_8ex}.}
\label{fig_normal_form}
\end{center}
\end{figure}

By Lemma~\ref{lem_linking}, the fundamental domain of 
$\U_{(p,q)}^m$ is as in Fig.~\ref{fig_normal_form}, where
we have cut along a spanning disc of $\beta$. 
To form 
$\U_{(p,q)}^m$, take $m$ copies of this fundamental domain
joined cyclically end-to-end: see, e.g., Fig.~\ref{fig_mainthma}.

\subsection{Subtemplates}

The (very complicated) normal form $\U^m_{(n_i,n_{i+1})}$ is tame 
compared to the full template $\U^m_{(n_1,\ldots,n_k)}$. Fortunately, 
we do not need to consider the most general case. 

\begin{figure}[hbt]
\begin{center}
\psfragscanon
\psfrag{a}[][]{}
\psfrag{b}[][]{}
\psfrag{c}[][]{}
\psfrag{z}[][]{\LARGE $\Rightarrow$}
\psfrag{p}[][]{\LARGE$I_2$}
\psfrag{q}[][]{\LARGE$I_1$}
\psfrag{r}[][]{\LARGE$J_3$}
\psfrag{s}[][]{\LARGE$J_2$}
\psfrag{t}[][]{\LARGE$J_1$}
\psfrag{P}[r][c]{\LARGE$g(I_2)$}
\psfrag{S}[l][c]{\LARGE$g(J_2)$}
\includegraphics[angle=0,width=5in]{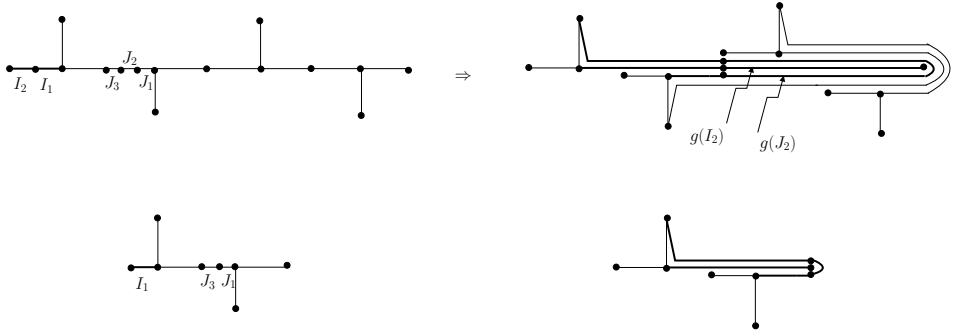}
\caption{$g = g_{(2,2,1)}$ [top]; 
$(S_1 \cup S_2)'$ [lower left]; 
$g((S_1 \cup S_2)')$, which equals $g_2(g_1(S_1 \cup S_2))$ [lower right].}
\label{fig_sub_lem}
\end{center}
\end{figure}

\begin{lem}
\label{lem_SubT}
Suppose $i$ is odd [resp. even]. Then $\T_{(n_i,n_{i+1})}$ [resp.
its mirror image] is a subtemplate of $\T_{(n_1,\ldots,n_k)}$.
\end{lem}
\pf
We prove this result on the level of the train track (by subdividing
and removing extraneous portions), and then suspend this to obtain
the result for subtemplates. We only prove the case $n_i, n_{i+1}>1$, 
for other cases are similar. 
Let $G=G_{(n_1, \cdots, n_k)}$ and $g=g_{(n_1, \cdots, n_k)}$. 
For the proof of the lemma, it is enough to show that 
there exist subdivisions 
$I_2I_1=e_i^{n_i}$, $J_3J_2J_1=e_{i+1}^{n_{i+1}}$ of $G$ 
that satisfy the following: 
\begin{quote}
Remove $I_2$ and $J_2$ from $S_i$ and $S_{i+1}$ respectively. 
Then 
$(S_i \cup S_{i+1})'$ denotes the graph obtained from 
$(S_i \cup S_{i+1})$ by identifying the 
right end point of $J_3$ and 
the left end point of $J_1$.
Then $g((S_i \cup S_{i+1})')=
g_{i+1}(g_{i}(S_i \cup S_{i+1}))$, see Fig.~\ref{fig_sub_lem}. 
\end{quote}
We demonstrate the case $i=1$, the proof for other cases 
being similar. Denote by $c_i$ the following edges of the train track: 
$c_i= e_i^0$ if $n_i = 1$ and $c_i = \overline{e_i^{n_i}} e_i^0$ 
if $n_i > 1$. 
By induction on $k$, one demonstrates that 
$g(e_1^{n_1})=e_1^0 c_2 \cdots c_k$, 
$g(e_2^{n_2})=e_2^0 c_3 \cdots c_{k-1} c_k \overline{c_k}~\overline{c_{k-1}} 
\cdots \overline{c_2} \overline{e_1^0}e_1^1$, and 
$g(e)=g_{2}(g_{1}(e))$ for $e \in S_1 \cup S_2 \setminus 
\{e_1^{n_1}, e_2^{n_2}\}$.  
Take subdivisions $I_2I_1=e_1^{n_1}$ and $J_3J_2J_1=e_2^{n_2}$ 
such that $g(I_1)=e_1^0 c_2$, $g(J_3) = e_2^0$, and 
$g(J_1)=\overline{c_2} \overline{e_1^0}e_1^1$. 
Since $g(I_1)=g_2(g_1(e_1^{n_1}))$, 
$g(J_3)g(J_1)= g_2(g_1(e_2^{n_2}))$ and 
$g(e)=g_{2}(g_{1}(e))$ for $e \in S_1 \cup S_2 \setminus 
\{e_1^{n_1}, e_2^{n_2}\}$, we have 
$g((S_1 \cup S_2)')=g_2(g_1(S_1 \cup S_2))$. 
\qed

\begin{lem}
\label{lem_SubU}
Suppose $i$ is odd [resp. even]. Then $\U^m_{(n_i,n_{i+1})}$ [resp.
its mirror image] is a subtemplate of $\U^m_{(n_1,\ldots,n_k)}$.
\end{lem}
\pf
We consider the case of $i$ odd and $n_i, n_{i+1}>1$; the 
other cases share identical proofs. 
The key step in this lemma is to determine the linking properties 
of $\beta$ with the subtemplate obtained from Lemma~\ref{lem_SubT}.

Note that $g_{(n_1, \cdots, n_k)}$ has three corresponding 
fixed points $q_{\l}', q_c', q_r'$, where $q_{\l}'$ and $q_r'$ are 
$(n_i+1)$ vertex and $(n_{i+1}+1)$ vertex of $S_i$ and $S_{i+1}$ 
respectively, and $q_c' \in e_{i+1}^{n_{i+1}}$.
Denote by $k_{\l}'$, $k_c'$, $k_r'$ the three unknots given 
by suspending these points. Let $\beta$ and $\beta'$ be the closures of 
$b_{(n_1,n_2)}$ and $b_{(n_1, \cdots, n_k)}$ respectively. 
As in the proof of Lemma~\ref{lem_linking},  
one shows that 
$\lk(k_{\l}', \beta')=\lk(k_{\l},\beta) = +1$, 
$\lk(k_{c}', \beta')=\lk(k_{c},\beta) = 0$, and 
$\lk(k_{r}', \beta')=\lk(k_{r}, \beta)= -1$, 
see Fig.~\ref{fig_linking_num4}(b). 
This together  with Lemma~\ref{lem_SubT} implies 
Lemma~\ref{lem_SubU}. 
\qed

%

\section{The Main Theorems}
\label{sec_maintheorems}

We now have assembled all the ingredients to prove the 
Main Theorems. 

\subsection{Main Theorem (c)}

We begin with the easiest and least general of the three Main Theorems.
A proof can be obtained implicitly by a combination of Lemma~\ref{lem_SubU}
with a geometric argument from \cite[pp. 444-5]{G}. However, to 
set a pattern for the proofs of Main Theorems (a) and (b), we
give a different proof of (c).

\begin{thm}
\label{thm_SimpleTrack}
Let $b=\braid$ where $n_j=n_{j+1}=1$
for some $j$. Then the closure of $b^m$ is a knot/link 
in $S^3$ with the universal fibration property for every $m\ge 2$.  
\end{thm}
\pf
 From Lemma~\ref{lem_SubU}, we may restrict attention to the 
subtemplate $\U^m_{(1,1)}$ within $\U^m_{(n_1,\cdots,n_k)}$. 
The template $\T_{(1,1)}$ is given in cartoon form in 
 Fig.~\ref{fig_8ex}, along with a representation of 
the fundamental domain of $\U^m_{(1,1)}$. 
Our strategy is to identify three unknots on $\U^m_{(1,1)}$ ---
$\kappa$, $\kappa'$, and $\kappa''$ --- to which we can apply
Theorem~\ref{thm_UT}. 

Specifically, using the symbolic representation of the orbits
as periodic words, let $\kappa=\word{\y}$, $\kappa'=\word{\x^m}$, 
and $\kappa''=\word{\z^m}$. Recall that $\word{\x^m}$ means to follow
the periodic orbit which traverses the $\x$-strip $m$ times. 
While there is an ambiguity in 
the symbolic representation (in terms of which fundamental 
domain the symbol is on), the only choice is for $\word{\y}$, 
and by cyclic symmetry, the result is independent of the choice. 

\begin{figure}[hbt]
\begin{center}
\psfragscanon
\psfrag{+}[][]{\Large $+$}
\psfrag{-}[][]{\Large $-$}
\includegraphics[angle=0,width=5.0in]{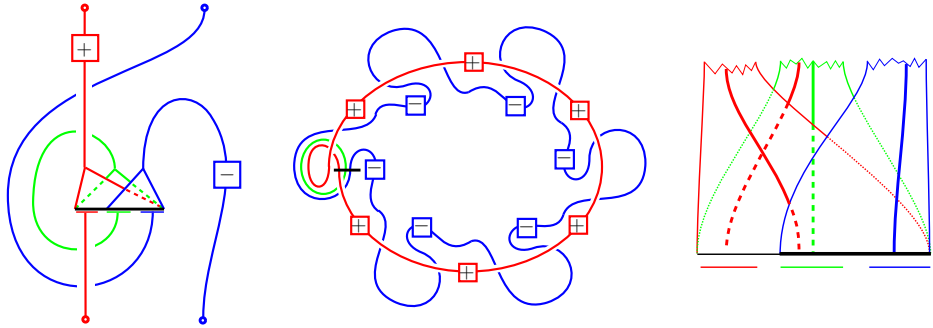}
\caption{A cartoon representation of the fundamental domain of 
$\U^m_{(1,1)}$ [left]; the unknots $\kappa$, $\kappa'$, and $\kappa''$
in the case $m=6$ [center]; all three orbits pass through the same 
branch line [right].} 
\label{fig_mainthma}
\end{center}
\end{figure}

It is trivial to see from Fig.~\ref{fig_mainthma}[left] that
(1) these three curves are separable unknots; and (2) $\kappa$ has
zero twist, $\kappa'$ has positive twist, and $\kappa''$ has
negative twist. The only other ingredient needed to invoke
Theorem~\ref{thm_UT} is that the three orbits meet at
a branchline as in Fig.~\ref{fig_three_orbits}. This is not the case,
since the strip $\z$ does not fully cover the branchline.
To remedy this, we replace $\word{\x^m}$ with $\word{\y\x^m}$.
This orbit is still a separable unknot (using 
\cite[Prop. 3.1.19]{GHS}) whose (already positive) twist has increased 
by two and which now passes through the triple branchline in 
precisely the manner of Fig.~\ref{fig_three_orbits}. 
Hence, by Theorem~\ref{thm_UT}, $\U^m_{(1,1)}$ is universal.
\qed

\subsection{Main Theorem (b)}

The following theorem will allow us to completely classify the 
fibrations of 2-bridge knot complements in the next section. 

\begin{thm}
\label{thm_TwoBridgeKnot}
The closure of $(\braid)^2$ has the universal fibration property.
\end{thm}
\pf 
To prove Theorem~\ref{thm_TwoBridgeKnot}, 
it suffices to show that 
$\U^2_{(p,q)}$ is universal for all $p$ and $q$. 
We assume $p,q >1$, the other cases being both simpler and 
similar to Main Theorem (c). The strategy, as in the previous 
theorem, is to find three well-chosen unknots $\kappa,\kappa',\kappa''$ 
on $\U^2_{(p,q)}$. 

\begin{figure}[hbt]
\begin{center}
\psfragscanon
\psfrag{a}[][]{\Large $+$}
\psfrag{c}[][]{\Large $-$}
\includegraphics[angle=0,width=5in]{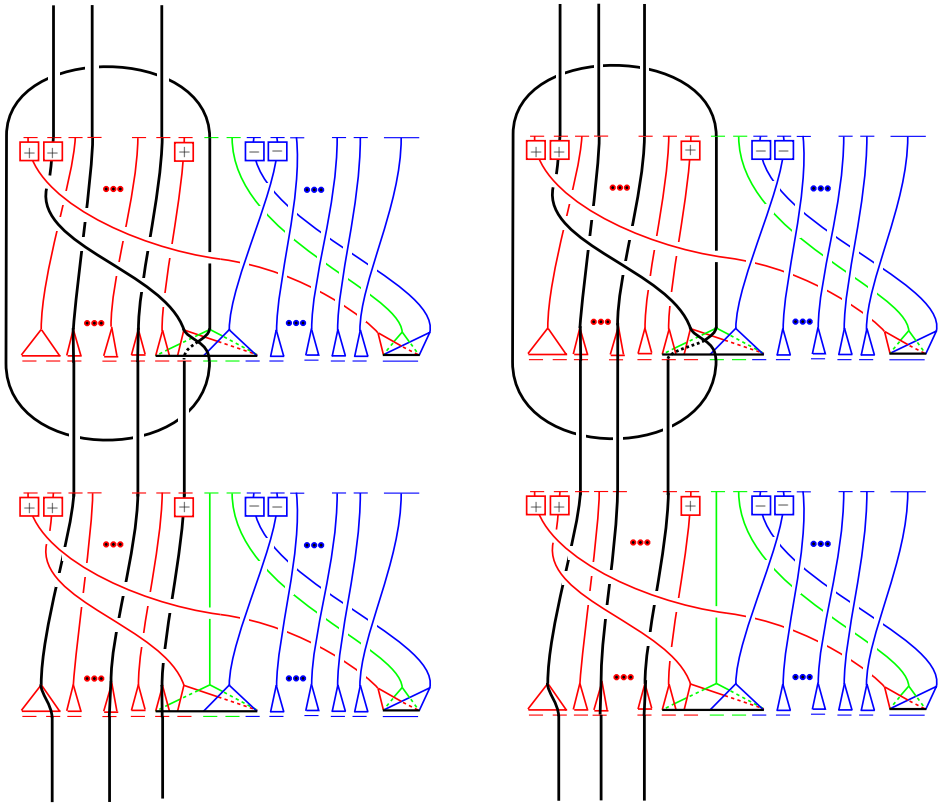}
\caption{The orbit $\kappa'$ in Main Theorem (b); 
$p$ is odd [left]; even [right].} 
\label{fig_positive}
\end{center}
\end{figure}

\begin{figure}[hbt]
\begin{center}
\psfragscanon
\psfrag{a}[][]{\Large $+$}
\psfrag{c}[][]{\Large $-$}
\includegraphics[angle=0,width=5in]{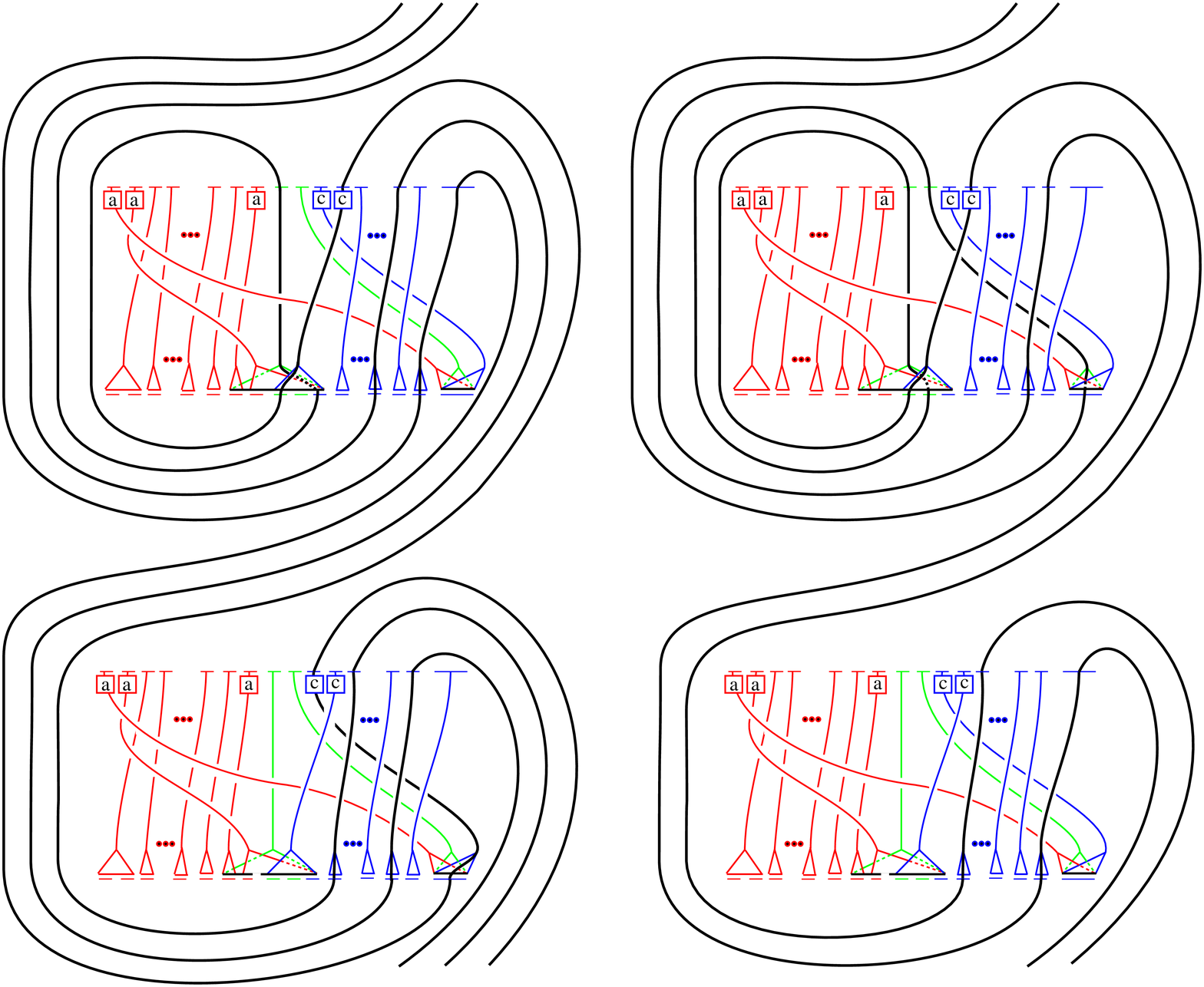}
\caption{The orbit $\kappa''$ in Main Theorem (b): 
$q$ is odd [left]; even [right].} 
\label{fig_negative}
\end{center}
\end{figure}

\begin{figure}[hbt]
\begin{center}
\psfragscanon
\psfrag{a}[b][]{\Huge $\kappa'$}
\psfrag{b}[b][]{\Huge $\kappa$}
\psfrag{c}[b][]{\Huge $\kappa''$}
\includegraphics[angle=0,width=4.75in]{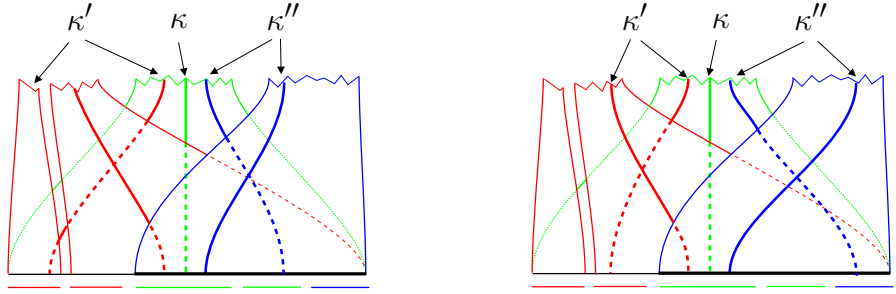}
\caption{A closeup of the branch line in the proof of Main Theorem (b);
$p$ and $q$ even [left]; $p$ and $q$ odd [right].} 
\label{fig_branch}
\end{center}
\end{figure}

As before, choose $\kappa:=(\y)^{\infty}$ the untwisted unknot 
separable from all other orbits on the template. 
For $p$ odd (resp. even), define
\[
\kappa':=\word{\y\x_0\x_1\x_2\cdots\x_{p-1}\x_p}
\quad
\left(\mbox{resp.}\  
\kappa':=\word{\y\x_1\x_2\cdots\x_{p-1}\x_p} 
\right), 
\]
see Fig.~\ref{fig_positive}. 
%
%
For $q$ odd (resp. even), define
\[
\kappa'':=\word{\y\z_{q}\z_0\z_1\cdots\z_{q-1}}
\quad
\left(\mbox{resp.}\  
\kappa'':=\word{\y\z_{q+1}\z_0\z_1\cdots\z_{q-1}}
\right), 
\]
see Fig.~\ref{fig_negative}. 
Note that there is an ambiguity as to which fundamental domain 
the orbit segment $y$ of $\kappa'$ and $\kappa''$ passes through. 
We specify as follows: take $\kappa'$ and $\kappa''$ such that 
the orbit segment $y$ passes through the fundamental domain where 
$\kappa$ lies. The figures imply that 
these are separable unlinked unknots with 
$\tau(\kappa) = 0$, $\tau(\kappa')>0$, and $\tau(\kappa'')<0$. 
These unknots $\kappa, \kappa'$ and $\kappa''$ pass through 
the triple branchline as in Fig.~\ref{fig_branch}. 
Hence, by Theorem~\ref{thm_UT}, $\U^2_{(p,q)}$ is universal.
\qed

\subsection{Main Theorem (a)}

For part (a) of the Main Theorem --- that the universal fibration 
property holds for $m$ sufficiently large --- we give a more precise
formulation.

\begin{thm}
Let $N:=\min_i(\max\{n_i,n_{i+1}\})$.  
The closure of $(\braid)^m$ has the 
universal fibration property for all $m\geq N(N-1)$. 
\end{thm}
\pf 
By Lemma~\ref{lem_SubU}, it is enough to show that 
$\U_{(p,q)}^m$ is universal for each $m\geq N(N-1)$, 
where $N=\max\{p,q\}$. As before, we assume $p$ and $q$ 
greater than $1$. 

Let $m=pt'+s'=qt''+s''$ be Euclidean integer factorizations
of $m$; hence $0\leq s'<p$ and $0\leq s''<q$. As $m\geq N(N-1)$, 
it follows that $t'\geq s'$ and $t''\geq s''$.
Consider the following three periodic orbits on $\U_{(p,q)}^m$:
\begin{equation}
\begin{array}{rcl}
  \kappa &:=& \word{\y} \\
  \kappa' &:=& 
    \word{(\y\x_1 \cdots \x_p)^{t'-s'}
    (\x_0\x_1 \cdots \x_p)^{s'}} \\
  \kappa'' &:=& 
    \word{(\z_{q+1}\z_0\z_1\cdots\z_{q-1})^{t''-s''}
    (\z_{q}\z_0 \z_1 \cdots \z_{q-1})^{s''}} 
\end{array} .
\end{equation}
The orbit $\kappa$ is, as before, an untwisted unknot which lies
within one fundamental domain of $\U_{(p,q)}^m$. The orbit 
$\kappa'$ passes through $p(t'-s')+(p+1)s'=pt'+s'=m$ fundamental
domains of $\U_{(p,q)}^m$. Except for the occasional excursions 
about the $\y$ orbit (which do not change the knot type), 
the orbit $\kappa'$ can be thought of as a 1-strand braid about
the braid axis of $\U_{(p,q)}^m$: it is therefore an unknot. 
All twisting in the $\x$-strips is of a positive nature, so 
that $\kappa'$ is an unknot with positive twist. 
A similar argument shows that $\kappa''$ passes through 
$q(t''-s'')+(q+1)s''=m$ fundamental domains and is an unknot
with negative twist separable from $\kappa'$. 

\begin{figure}[hbt]
\begin{center}
\psfragscanon
\psfrag{a}[b][]{\Large $\kappa'$}
\psfrag{b}[b][]{\Large $\kappa$}
\psfrag{c}[b][]{\Large $\kappa''$}
\includegraphics[angle=0,width=2.95in]{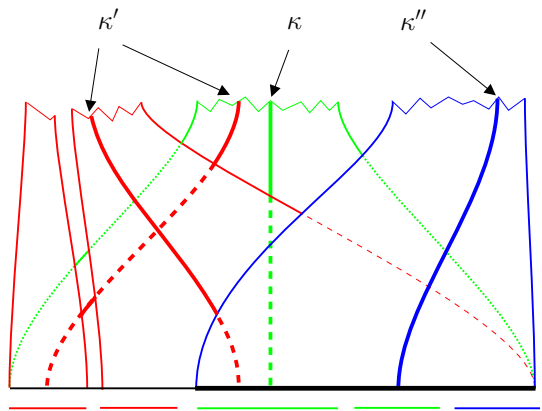}
\caption{A closeup of the branch line in the proof of Main Theorem (a).}
\label{fig_mainthmc}
\end{center}
\end{figure}

These orbits are a triple of separable unknots which satisfy the 
conditions  of Theorem~\ref{thm_UT}, as illustrated in 
Fig.~\ref{fig_mainthmc}. Hence $\U_{(p,q)}^m$ is universal. 
\qed

%

\section{Concluding remarks}
\label{sec_conc}

We have displayed a large family of fibred knots/links 
whose fibrations force all knot and link types in the complement.
It remains an open question how large this family can be:

\begin{question}
\label{q_universal}
Is every fibred knot/link of pseduo-Anosov type universally fibred ?
\end{question}

We conclude this paper with a few remarks on this question.

\subsection{2-bridge knots}

Recall that a \df{2-bridge} knot is one which can be isotoped
so as to have precisely four critical points (two minima and 
two maxima) with respect to the standard height function in $\real^3$.
We can affirmatively answer Question~\ref{q_universal} in this
category.

\begin{figure}[hbt]
\begin{center}
\includegraphics[angle=0,width=5in]{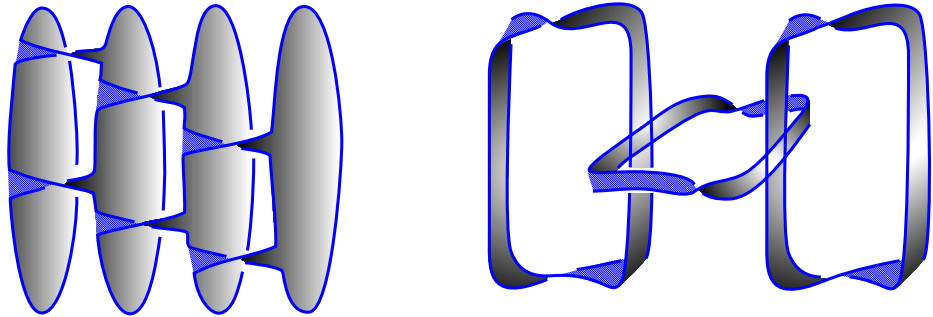}
\caption{Plumbings of Hopf bands [right] can be arranged as 
a closed braid of the form $(b_{n_1, \cdots, n_k})^2$ [left].}
\label{fig_2Bridge}
\end{center}
\end{figure}

\begin{cor}
\label{cor_2bridge}
Every fibred non-torus 2-bridge knot 
has the universal fibration property.
\end{cor}
\pf
The work of Gabai and Kazez \cite{GK} explicitly shows that the 
class of fibred 2-bridge knots is precisely that class of knots whose
fibres are a chain of plumbed Hopf bands, and that the monodromy 
is of pseudo-Anosov type if and only if the Hopf bands are not 
all of the same sign, which is to say, if and only if it is not a
torus knot.
We claim that such a chain of plumbed Hopf bands has a closed braid form 
of $(\braid)^2$. The proof is inductive on the number of 
Hopf bands and illustrated in Fig.~\ref{fig_2Bridge}. The key step is 
noting that the fibre for the closure of $(\braid)^2$ is a stack of
discs attached by a pairs of half-twisted bands: this is precisely
a chain of plumbed Hopf bands. From Theorem~\ref{thm_TwoBridgeKnot}, we 
have that this knot must be universally fibred. 
\qed

\subsection{Connected sums}

We can enlarge the class of universally fibred knots considerably 
by taking connected sums.

\begin{prop}
\label{prop_ConnectSum}
The connected sums of any fibred knot/link with any universally fibred
knot/link is universally fibred. 
\end{prop}

\pf
In \cite[Corollary 1.4.]{Gabai2}, it is shown that the Murasugi sum 
of fibre surfaces is ``natural''. Specifically, if 
$K_i$ is a fibred knot/link ($i=1,2$) with fibre surface 
$\Sigma_i$ and monodromy $\Phi_i$ satisfying 
$\Phi_i|_{\partial{\Sigma_i}}=id$, 
then the Murasugi sum of $\Sigma_1$ and $\Sigma_2$ along a $2k$-gon $D$ 
common to $\Sigma_1$ and $\Sigma_2$ is fibre surface and 
has monodromy $\tilde\Phi_1\circ\tilde\Phi_2$, where $\tilde\Phi_i$ is 
equal to $\Phi_i$ on $\Sigma_i$ and the identity map elsewhere. 

Even more importantly, the proof of this result uses vector
fields transverse to the fibration to achieve the result. Thus, 
we actually have that a transverse vector field of the 
Murasugi sum decomposes naturally. 

The simplest kind of Murasugi sum (where the attachment is 
along a 2-gon in the spanning surfaces) is precisely the
connected sum. Choose the monodromy $\Phi_1$ of the 
universally fibred knot to be the pseudo-Anosov representative,
isotoped so as to be the identity map on the 
the boundary. Likewise, choose the monodromy $\Phi_2$ to be 
the Thurston canonical form isotoped to be the identity on the
boundary. The composition monodromy on the connected 
(in this case, Murasugi) sum is easily isotoped to have an invariant
curve separating the universal pseudo-Anosov component
from the Thurston form of $\Phi_2$. Thus, the theorem of 
Asimov and Franks implies that any transverse vector field must have
at least this set of periodic orbits --- all knots and links. 
\qed

As a corollary of the proof, we see that in certain instances, 
more general Murasugi sums (e.g., plumbing) can preserve the
universal fibration property as well. One would need to control
(1) the dynamics on the attaching region, to preserve the 
pseudo-Anosov property of the monodromy, and (2) the location of
the attaching region, so that it does not interfere with the 
universal subtemplate. 

\subsection{Too-twisted monodromies}


We believe that there are pseudo-Anosov fibred knots or links whose
monodromies are too ``twisted'' to admit a universal template
in the suspension. Our proposed example stems from the fact
that the proof of Main Theorem (c) requires a sufficiently
high order branched cover.

\begin{conj}
For $p$ and $q$ sufficiently large, the closure of 
$(b_{(p,q)})^3$ is not universally fibred.
\end{conj}

Our reason for 
conjecturing the non-universality stems from the difficulty
in finding enough unknotted orbits on which to build a 
universal subtemplate.

\end{document}